\documentclass{elsart_mod}

\usepackage[dvips]{graphicx}
\usepackage{amsmath}
\usepackage{amssymb}
\usepackage{enumerate}
\usepackage{natbib}
\usepackage{url}

\newcommand{\eref}[1]{(\ref{#1})}

\newcommand{\E}{\textrm{E}}
\newcommand{\argmax}{\operatornamewithlimits{argmax}}
\newcommand{\argmin}{\operatornamewithlimits{argmin}}

\begin{document}
\begin{frontmatter}
\title{Approximate cost-efficient sequential designs for binary response models with application to switching measurements}
\author{Juha Karvanen}
\address{
National Public Health Institute,\\
Mannerheimintie 166, 00300 Helsinki, Finland}
\ead{juha.karvanen@ktl.fi}

\begin{abstract}
The efficiency of an experimental design is ultimately measured in terms of time and resources needed for the experiment. Optimal sequential (multi-stage) design is studied in the situation where each stage involves a fixed cost.
The problem is motivated by switching measurements on superconducting Josephson junctions. In this quantum mechanical experiment, the sequences of current pulses are applied to the Josephson junction sample and a binary response indicating the presence or the absence of a voltage response is measured. The binary response can be modeled by a generalized linear model with the complementary log-log link function. The other models considered  are the logit model and the probit model. For these three models, the approximately optimal sample size for the next stage as a function of the current Fisher information and the stage cost is determined. The cost-efficiency of the proposed design is demonstrated in simulations based on real data from switching measurements. The results can be directly applied to switching measurements and they may lead to substantial savings in the time needed for the experiment.
\end{abstract}
\end{frontmatter}

\textsf{J. Karvanen, Approximate cost-efficient sequential designs for binary response models with application to switching measurements. Computational Statistics \& Data Analysis, doi:10.1016/j.csda.2008.10.018, 2008.\\
\\
The older title of the paper was ``Cost-efficient sequential designs for binary response models''.}

\section{Introduction}
In binary response experiments carried out in various fields of experimental science, the researcher chooses the value of a covariate variable and measures a binary response. The measurement data is modeled by a binary response model where the probability of success or failure is assumed to be a monotone function of the covariate. Three parametric models, the logit, the probit and the complementary log-log (cloglog) model are frequently used to model the dependence between the binary response $Y$ and the continuous covariate $x$. All these models can be presented in the framework of generalized linear models \citep{McCullagh:glm}
\begin{equation} \label{eq:binary}
P(Y=1)=\textrm{E}(Y)=F(ax+b),
\end{equation}
where the response curve $F$ is a cumulative distribution function (cdf) and the $a$ and $b$ are the parameters of the model to be estimated. The three response curves commonly used are:
logistic distribution for the logit model
\begin{equation} \label{eq:logit}
F(ax+b)=\frac{\exp(ax+b)}{1+\exp(ax+b)},
\end{equation}
normal distribution for the probit model
\begin{equation} \label{eq:probit}
F(ax+b)=\Phi(ax+b),
\end{equation}
where $\Phi$ is the cdf of the standard normal distribution,
and the Gompertz distribution for the cloglog model
\begin{equation} \label{eq:cloglog}
F(ax+b)=1-\exp(-\exp(ax+b)).
\end{equation}

The problem of optimal design arises when the covariate values need to be chosen in optimal manner to estimate some parameters of interest. In binary response models~\eref{eq:binary}, the parameters of interest are $a$ and $b$ that determine the slope and the location of the response curve. Several optimality criteria have been proposed in the statistics literature. Among the most popular are D-optimal designs, which minimize the generalized variance. The minimization of the generalized variance equals to the maximization of the determinant of the Fisher information. The expected Fisher information of the data $x_{i}$, $i=1,\ldots,n$ under the binary response model~\eref{eq:binary} is defined as
\begin{equation} \label{eq:Fisherinf}
\mathbf{J}=
\begin{pmatrix}
\displaystyle \sum_{i=1}^{n}g(z_{i})x_{i}^{2} & \displaystyle \sum_{i=1}^{n}g(z_{i})x_{i}\\
\displaystyle \sum_{i=1}^{n}g(z_{i})x_{i} & \displaystyle \sum_{i=1}^{n}g(z_{i})
\end{pmatrix},
\end{equation}
where $z_{i}=ax_{i}+b$ and
\begin{align}
g(z)=&\frac{e^{-z}}{(1+e^{-z})^{2}} \quad \textrm{for logit model~\eref{eq:logit}},\\
g(z)=&\frac{(\Phi^{'}(z))^{2}}{\Phi(z)(1-\Phi(z))} \quad \textrm{for probit model~\eref{eq:probit} and}\\
g(z)=&\frac{e^{2z}}{e^{\exp(z)}-1} \quad \textrm{for cloglog model~\eref{eq:cloglog}}.
\end{align}
Now the D-criterion can be defined as the square root of the determinant of Fisher information \eref{eq:Fisherinf}
\begin{equation}
D=\sqrt{\det(\mathbf{J})}.
\end{equation}
The (locally) D-optimal design is derived by maximizing $D$ under the assumptions that $n \rightarrow \infty$ and the true values of the the parameters $a$ and $b$ are known.  The derivation is non-trivial and can be found in the general form in \citep{Ford:canonicaloptimal}. The widely used logit model~\eref{eq:logit} is studied also by \citet{Abdelbasit:experimentaldesign,Minkin:optimal,Sitter:optimal,Mathew:optimal}.  The above references show that the D-optimal design is a two-point design for the logit, the probit and the cloglog models and the optimal covariate values $x^{*}$ can be solved from the equation $ax^{*}+b=z^{*}$, where the values of the canonical parameter $z^{*}$ are reported in Table~\ref{tab:Doptims}. Extensions to the models with more than one covariate are considered by \citet{Woods:designs}, \citet{Dror:sequential} and \citet{DortaGuerra:twolevel}.

\begin{table}[htb]
\caption{D-optimal covariate values for the models~\eref{eq:logit}, \eref{eq:probit} and \eref{eq:cloglog}. The optimal covariate values $x_{1}^{*}$ and $x_{2}^{*}$ are solved from equations $ax_{1}^{*}+b=z_{1}^{*}$ and $ax_{2}^{*}+b=z_{2}^{*}$. Columns $F(z_{1}^{*})$ and $F(z_{2}^{*})$ report the cdf values related to the optimal covariates, i.e. probabilities of response~1. \label{tab:Doptims} }
\begin{center}
\begin{tabular}{rcccc}
model & \multicolumn{4}{c}{Optimal covariates}\\
 & $z_{1}^{*}$ & $z_{2}^{*}$ & $F(z_{1}^{*})$ & $F(z_{2}^{*})$\\ \hline
logit  & -1.543 & 1.543 & 0.176 & 0.824\\
probit  & -1.138 & 1.138 & 0.128 & 0.872\\
cloglog  & -1.338 & 0.980 & 0.231 & 0.930\\
\end{tabular}
\end{center}
\end{table}

The efficient use of resources is not measured by the D-criterion but in terms of time and resources used for the experiment. Optimal sequential (multi-stage) design  is studied in the situation where there is a fixed cost for each stage. As a motivating example, switching measurements on superconducting Josephson junctions \citep{Josephson,single,pekola} are considered. In this quantum mechanical experiment, the experimenter varies the height of the applied current pulse and observes the presence or the absence of the voltage response generated by macroscopic quantum tunneling~\citep{legget}. An experimental design for switching measurements was proposed by \citet{optdesign}. The following properties are characteristic for switching measurements:
\begin{itemize}
\item The data acquisition is fast and the total number of observations can be high. In the experiment reported in~\citep{optdesign}, the total number of observations was over 100000 and the experiment took about 8 minutes. This means that there is no time for expert analysis during the experiment but automated processing must be used.
\item Sequential/multi-stage designs are applicable. There is, however, a cost related to the number of the stages which usually makes the fully sequential approach inefficient.
\item The measurement system remains stable only for a limited time because very low temperatures are needed for superconductivity. This emphasizes the need for cost-efficient designs.
\item Prior information on the model parameters may be very poor.
\item Differently from dose-response trials~\citep{Biedermann:restricteddesignspaces}, ethical considerations do not restrict the choice of the design.
\end{itemize}

Sensitivity to the modeling assumptions and the need of initial estimates are the well-known disadvantages of D-optimality. In switching measurements, however, there are good reasons to use D-optimality: First, it can be shown \citep{optdesign} that the cloglog model~\eref{eq:cloglog} provides a very good approximation for the probability of switching as a function of the height of the applied current pulse. In other words, in switching measurements, the true model is actually known. Second, the experiment is carried out sequentially, which means that poor initial estimates will not spoil the whole experiment. In sequential switching measurements, a number of current pulses with a certain height are generated, responses are measured for each pulse and finally the height of the pulse for the next stage is calculated. Because the configuration of the pulse generator takes time, it is not practical to change the height of the pulse after every pulse. The total cost equals to the time used for the experiment and consist of the time used for measuring and computing and the time used for the configuration of the pulse generator between the stages.

The optimal number of measurements per stage in sequential designs with a fixed stage cost is determined for the logit, the probit the cloglog models. The proposed solution is approximal and based on numerical analysis. It is assumed that the initial maximum likelihood estimates for the parameters are available. If the initial estimates are not available, they can be found, for instance, by the binary search algorithm \citep{inidesign}. After the initial estimation, the experiment continues as a sequential design where the D-optimality criterion is used. The optimal stage size, i.e., the optimal number of measurements per stage, is derived as a function of the stage cost $C_{S}$ and the observed value of $D$. The derived cost-efficient designs can be presented in simple functional form and are easily applied to practical experiments.

Most of the existing works on sequential design do not consider the cost of the experiment as a design criterion.
\citet{Wu:fullysequential} considered fully sequential designs where the measurements are done one by one. He criticized the earlier approaches: the stochastic approximation method \citep{RobbinsMonro} and the up-and-down method \citep{Dixon:upanddown} and proposed a sequential approach to the estimation of 100$p$th percentile. In this approach, the next covariate point $x_{n+1}$ is chosen so that $\hat{F}_{n}(x_{n+1})=p$, where $\hat{F}_{n}$ is the estimated cdf of logistic distribution. \citet{Ivanova:improvedupanddown} considered several designs from the family of up-and-down rules for the sequential allocation of dose levels to subjects in a dose-response study. The use of group up-and-down designs \citep{Gezmu:groupupanddown} implies that there is a cost related to the number of stages also in dose-finding experiments although this cost is not usually explicitly considered. \citet{McLeish:sequentialbioassay} studied sequential design in bioassay and used a cost function that included unit cost plus extra cost for each positive response. \citet{Sitter:optimaltwostage} and \citet{Sitter:twostagequantal} studied two-stage designs for the logit and the probit model and used the second stage to balance the first stage. The second stage design was found by numerical optimization. The results can be extended to the multi-stage designs. Theoretical results on asymptotics of sequential designs were provided in \citep{Chaudhuri:nonlinearexperiments} and \citep{Chaudhuri:efficientdesigning}. A review of Bayesian approach to sequential designs was given by \citet{Chaloner:Bayesianexperimental}. A recent major contribution to the sequential design was the procedure proposed by \citet{Dror:sequential}. Their procedure is based on a Bayesian analysis that exploits a discretization of the parameter space to efficiently represent the posterior distribution. \citet{Tekle:maximin} studied maximin D-optimal designs for binary longitudinal responses. \citet{twostage} considered two-stage designs for gene-disease association studies where the motivation for the design arises from cost considerations.

In Section~\ref{sec:costeffi}, cost-efficient sequential designs are derived assuming that the initial maximum likelihood estimates of the parameters exist. In Section~\ref{sec:application}, application to switching measurements is presented. Discussion in Section~\ref{sec:discussion} concludes the paper.

\section{Cost-efficient sequential designs} \label{sec:costeffi}

\subsection{Problem definition} \label{sec:problemdefinition}
Consider a sequentially performed experiment where the total cost of the experiment consists of two components: a cost related to the number of measurements and a cost related to the number of stages. Without loss of generality, the marginal cost of making one additional measurement is fixed as unity and the cost of having one additional stage is marked by $C_{S}$. Thus, the total cost of an experiment with total of $n$ measurements in $K$ stages would be $C=n+KC_{S}$. The number of measurements at stage $k$ is marked by $n_{k}$, $k=1,\ldots,K$. It follows that $\sum_{k=1}^{K}n_{k}=n$. Assume that after $k-1$ stages the maximum likelihood estimates $\hat{a}_{k-1}$ and $\hat{b}_{k-1}$ and the related Fisher information matrix $\mathbf{J}_{k-1}$ are obtained. On the basis of the estimates $\hat{a}_{k-1}$ and $\hat{b}_{k-1}$ and the results in Table~\ref{tab:Doptims}, the estimated D-optimal covariate points $x_{k1}^{*}$ and $x_{k2}^{*}$ are calculated.  The problem is to decide the optimal number of measurements $n_{k}$ to be made in the points $x_{k1}^{*}$ and $x_{k2}^{*}$. For practical reasons, it is assumed that $n_{k}$ is an even number so that the measurements can be divided equally between $x_{k1}^{*}$ and $x_{k2}^{*}$. Because $\hat{a}_{k-1}$ and $\hat{b}_{k-1}$ are only estimates of the true parameters, the covariate points $x_{k1}^{*}$ and $x_{k2}^{*}$ may be far from the actual truly D-optimal covariate points. If $n_{k}$ is small, $\hat{a}$ and $\hat{b}$ are updated quickly and the number of measurements needed is small but the number of stages needed is high. On the other hand if $n_{k}$ is large, the number of measurements needed is higher but the number of stages needed is smaller. Obviously, the optimal $n_{k}$ depends on the stage cost $C_{S}$ and the Fisher information $\mathbf{J}_{k-1}$.

The scheme described above approximates the cost structure of switching measurements. In switching measurements, there are actually three types of costs involved:
\begin{description}
\item[marginal unit cost] time needed for making one additional measurement with the present pulse height,
\item[configuration cost] time needed for changing the pulse height,
\item[computing cost] time needed for updating parameter estimates.
\end{description}
The marginal unit cost and the configuration cost remain almost constant during the experiment but the computing cost depends on the details of numerical optimization and may vary from stage to stage. Additional complications follow from the fact that the D-optimal design is a two-point design implying that the signal generator needs to be configured when moving from one design point to another. This does not prevent using the described scheme because the stage cost may be defined as a sum of the computing cost plus two times the configuration cost.

The problem of cost-efficient design may be expressed as an optimization problem with respect to a precision target or a budget constraint. If a precision target $\psi_{\textrm{target}}$ is used, the optimal $n_{k}$ is found as a solution to the following minimization problem
\begin{equation} \label{eq:analytical1}
\argmin_{n_{k}}\; \E\left(C(\mathbf{x},\mathbf{y})\,|\,\psi(\mathbf{x},\mathbf{y}) \geq \psi_{\textrm{target}}, \mathbf{x}_{k-1},\mathbf{y}_{k-1},n_{k}\right),
\end{equation}
where  $\mathbf{x}_{k-1}$ and $\mathbf{y}_{k-1}$ present all measurement data available after $k-1$ stages,
$C(\mathbf{x},\mathbf{y})$ is the cost of collecting measurement data $(\mathbf{x},\mathbf{y})$ and $\psi(\mathbf{x},\mathbf{y})$
is the value of the precision criterion calculated from the measurement data $(\mathbf{x},\mathbf{y})$. For instance, $D$ may be used as the precision criterion.
The expectations in minimization~\eref{eq:analytical1} are integrals over
all possible measurement data $(\mathbf{X}_{k},\mathbf{Y}_{k})$ available after $k$ stages. For data $(\mathbf{X}_{k}=\mathbf{x}_{k},\mathbf{Y}_{k}=\mathbf{y}_{k})$ the expected cost to reach the precision target is given by
\begin{equation}
k C_{S}+\sum_{i=1}^{k} n_{i}+\min_{n_{k+1},n_{k+2},\ldots} \E\left(C(\mathbf{x},\mathbf{y}) \,|\,\psi(\mathbf{x},\mathbf{y}) \geq \psi_{\textrm{target}}, \mathbf{x}_{k},\mathbf{y}_{k}\right),
\end{equation}
where the last term in the summation is the expected cost to reach the precision target on the condition of data $(\mathbf{x}_{k},\mathbf{y}_{k})$
when the optimal strategy is used for the rest of the experiment. In the other words, minimization~\eref{eq:analytical1} involves recursive calculations where the same minimization is performed for the subsequent stage sizes $n_{k+1},n_{k+2}, \ldots$.

Alternatively, if a budget constraint  $C_{\textrm{budget}}$ is used, the optimal $n_{k}$ is found as a solution to the following maximization problem
\begin{equation} \label{eq:analytical2}
\argmax_{n_{k}}\; \E\left(\psi(\mathbf{x},\mathbf{y})\,|\,C(\mathbf{x},\mathbf{y}) \leq C_{\textrm{budget}}, \mathbf{x}_{k-1},\mathbf{y}_{k-1},n_{k}\right).
\end{equation}
Maximization \eref{eq:analytical2} requires similar recursive calculations as minimization~\eref{eq:analytical1}.

In stochastic optimization \citep{Birge:stochasticprogramming}, the problems \eref{eq:analytical1} and \eref{eq:analytical2} may be characterized as multi-stage recourse problems where the number of the stages is not fixed. Because of the nonlinearity of the objective functions and the recursive structure of the problem, it seems difficult to find optimal $n_{k}$ in analytical means. The high number of branches in the recursion implies that straightforward simulations cannot solve the problem either, at least not in reasonable time. To overcome these difficulties, our solution uses numerical computation and simplifying approximations.

\subsection{Design paths in the $(D,C)$-plane}
\begin{figure}
  \includegraphics[width=\columnwidth]{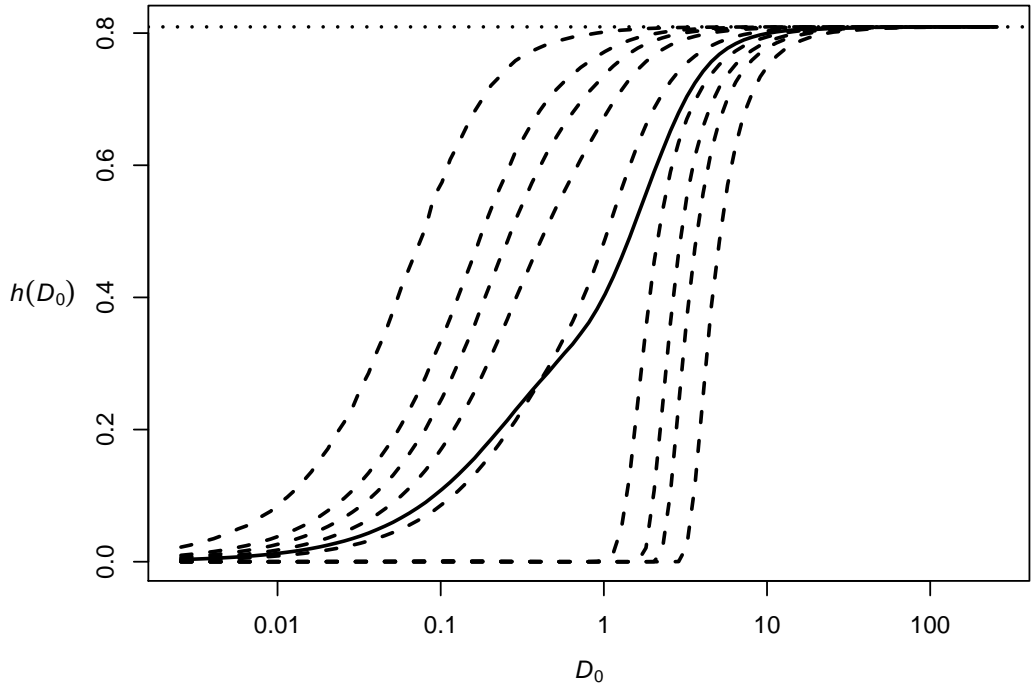}
  \caption{The expected information increase as a function of current information for the cloglog model when two additional measurements are made. The solid line represents the average of $h(D_{0})$ in 400000 simulation runs and the dashed lines represents 1st, 5th, 10th, 20th, 50th, 80th, 90th, 95th and 99th percentiles. The theoretical maximum of $h(D_{0})$ corresponding to the D-optimal two-point design is marked by the dotted line.}\label{fig:deltadcloglog}
\end{figure}

Because it is difficult to derive cost-efficient sequential designs analytically, approximations and simulations are used to obtain designs that are approximately cost-efficient. The first approximation made is to compress the information matrix $\mathbf{J}$ into the square root of its determinant $D$. This allows to compare candidate designs in the $(D,C)$-plane. Each design is characterized by the number of measurements $n_{1},n_{2},\ldots,n_{K}$ per stage and can be drawn as a path in the $(D,C)$-plane so that each point of the path represents the expected value of $D$ as a function of $C$. At the beginning of each stage there is a vertical jump in the path corresponding to the stage cost $C_{S}$. After that $D$ is expected to increase with a constant slope as $C$ increases. If two additional measurements are made, the expected change in $D$ is denoted $h(D_{0})$, where $D_{0}$ is the value of $D$ at the beginning of the stage. The number of the additional measurements is two because the D-optimal designs are two-point designs. An approximation for the function $h(D_{0})$ is obtained by means of simulation. The idea is to simulate the distribution of $\hat{a}$ and $\hat{b}$ on the condition that the expected square root of the determinant of the information matrix is $D_{0}$. Without loss generality it is fixed $a=1$ and $b=0$.  The information matrix is taken to be
\begin{equation}
\mathbf{J}_{0}=\frac{D_{0}}{D^{*}}\mathbf{J}^{*},
\end{equation}
where
\begin{equation}
\mathbf{J}^{*}=
\begin{pmatrix}
g(z_{1}^{*})(z_{1}^{*})^{2}+g(z_{2}^{*})(z_{2}^{*})^{2} & g(z_{1}^{*})z_{1}^{*}+g(z_{2}^{*})z_{2}^{*}\\
 g(z_{1}^{*})z_{1}^{*}+g(z_{2}^{*})z_{2}^{*} & g(z_{1}^{*})+g(z_{2}^{*})
\end{pmatrix},
\end{equation}
which is the information matrix of the D-optimal covariates, and $D^{*}=\sqrt{\det{J^{*}}} \approx 0.8094$. An approximation that $\hat{a}$ and $\hat{b}$ are normally distributed with the mean vector $(a,b)$ and the covariance matrix $\mathbf{J}_{0}^{-1}$ is made. In the other words, it is assumed that the estimates $\hat{a}$ and $\hat{b}$ originate from previous D-optimal measurements. The value of $\det(\mathbf{J}(\hat{a},\hat{b}))$ is calculated for each generated pair of $\hat{a}$ and $\hat{b}$
\begin{equation}
\det(\mathbf{J}(\hat{a},\hat{b}))=\sqrt{g(z_{1})g(z_{2})(z_{1}-z_{2})^{2}},
\end{equation}
where
\begin{equation}
z_{1}=\frac{z_{1}^{*}-\hat{b}}{\hat{a}} \textrm{ and } z_{2}=\frac{z_{2}^{*}-\hat{b}}{\hat{a}}.
\end{equation}
The expected change $h(D_{0})$ is then obtained as average over the values of $\det(\mathbf{J}(\hat{a},\hat{b}))$. All numerical calculations are done using R \citep{R}. The obtained function $h(D_{0})$ is plotted in logarithmic scale in Figure~\ref{fig:deltadcloglog} for the cloglog model. When $D_{0}$ increases, $h(D_{0})$ approaches its theoretical maximum $0.8094$ that corresponds to making the two additional measurements at the D-optimal points $z_{1}^{*}$ and $z_{2}^{*}$. The graphs for the logit and the probit model are similar (not shown).

In order to use the expected change of the D-criterion in the optimization, the function $h(D_{0})$ need to be interpolated. By looking Figure~\ref{fig:deltadcloglog} it is found that $h(D_{0})$ can be roughly approximated by logistic function $D^{*}/(1+\exp(\eta+\theta \log(D_{0})))$, where $\eta$ and $\theta$ are parameters estimated from the simulated data. Using this approximation as the basis of the model, a generalized additive model \citep{Hastie:gambook,Hastie:gamRcode} is fitted to the simulated data. In the model, smoothing splines with the degree of freedom 8 are used to model the nonlinearities remaining after the logistic approximation. The estimated model fits very well to the simulated data so that in Figure~\ref{fig:deltadcloglog} it would not be easy to separate the observed $h(D_{0})$ from the modeled $h(D_{0})$ (the maximum absolute difference between the observed and modeled values of $h(D_{0})$ is 0.016).

\subsection{Comparing design paths} \label{sec:comparingpaths}
\begin{figure}
  \includegraphics[width=\columnwidth]{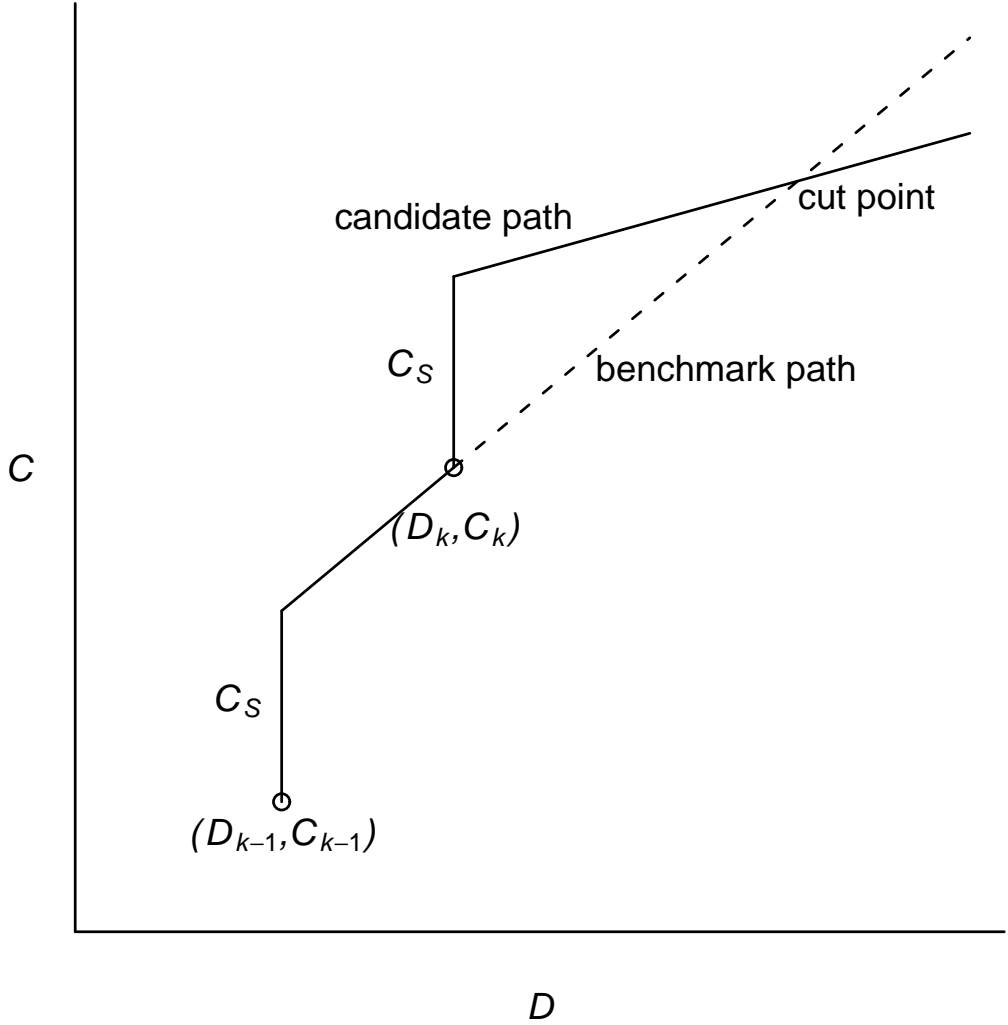}
  \caption{Schematic illustration of comparison of design paths. The benchmark path marked by dashed line starting from $(D_{k-1},C_{k-1})$ corresponds to completing the experiment in one stage whereas the candidate path path (solid line) corresponding to having first stage of $n_{k}$ measurements and then completing the experiment in one stage. The candidate design has an additional cost of $C_{S}$ compared with the benchmark design but because of the additional information obtained at the first stage, the increase in the information is faster at the second stage. The cut point of the paths is given by~\eref{eq:cutpoint} and the $D$-coordinate of the cut point serves as a criterion for comparing candidate paths with different $n_{k}$.}\label{fig:idea}
\end{figure}

Next the optimal number of measurements $n_{k}$ for the next stage are calculated providing that the current point in the $(D,C)$-plane is $(D_{k-1},C_{k-1})$ and the stage cost is $C_{S}$. For the calculation a two-stage approximation that evaluates a large number of candidate paths after two stages is used. The candidate paths are chosen such a way that the potential values of $n_{k}$ are covered.  It is assumed that our precision target (defined in terms of $D$) and budget restriction (defined in terms of $C$) do not restrict the evaluation, i.e. the target and the restriction are not met after two additional stages.  The idea of the calculation is illustrated in Figure~\ref{fig:idea}. The benchmark path (dashed line) is a line that goes through the point $(D_{k-1},C_{k-1}+C_{S})$ and has slope $2/h(D_{k-1})$. The benchmark path corresponds to the decision completing the experiment after having only one additional stage. The candidate paths have two stages. After the first stage the expected location is $(D_{k-1}+n_{k} h(D_{k-1})/2,C_{k-1}+C_{S}+n_{k})$. For the second stage the slope is  $2/h(D_{k})$. Because $h(D_{k})>h(D_{k-1})$, the candidate paths cut the benchmark path at point
\begin{align} \label{eq:cutpoint}
D=&D_{k-1}+\frac{n_{k} h(D_{k-1})}{2}+\frac{C_{S}h(D_{k-1})h(D_{k})}{2(h(D_{k})-h(D_{k-1}))},\\ \nonumber
C=&C_{k-1}+n_{k}+\left(1+\frac{h(D_{k})}{h(D_{k})-h(D_{k-1})} \right)C_{S}.
\end{align}
The cut points are calculated for a large number of candidate paths. The candidate path that has the earliest cut-point (the smallest $D$) with the benchmark path is chosen as the optimal one. The corresponding number of measurements $n_{k}$ is taken as optimal for $D_{k-1}$ when the stage cost is $C_{S}$.

\subsection{Approximately optimal sequential designs}
\begin{figure}
  \includegraphics[width=\columnwidth]{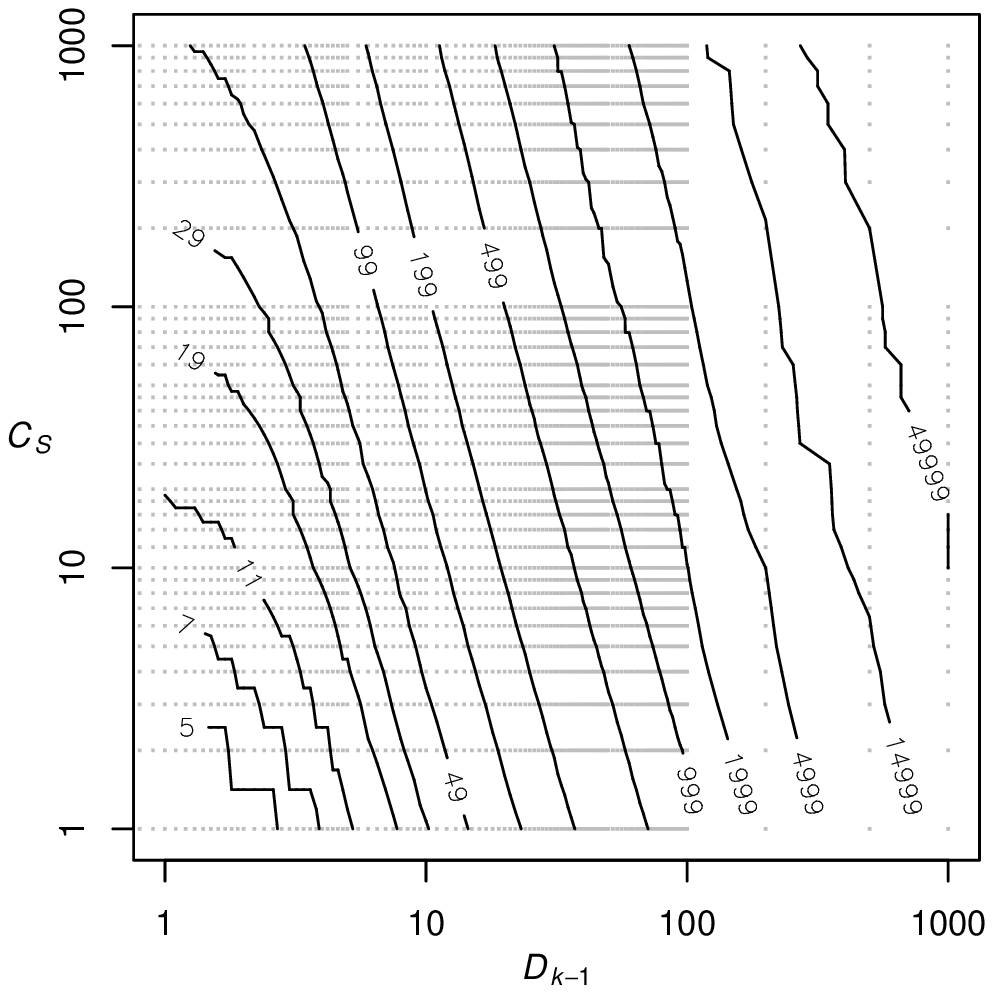}
  \caption{Optimal number of measurements $n_{k}$ for the next stage as a function of current information $D_{k-1}$ and stage cost $C_{S}$ for the cloglog model. The dots indicate the 3640 $(D_{k-1},C_{S})$ pairs for which the value of $n_{k}$ was determined. The contours are calculated from these numbers.} \label{fig:contourcloglog}
\end{figure}

In order to find a general model for the optimal $n_{k}$ as a function of the current information $D_{k-1}$ and the stage cost $C_{S}$, the optimal $n_{k}$ is calculated for a large number of $(D_{k-1},C_{S})$ pairs using the two-stage approximation described in Section~\ref{sec:comparingpaths}. The current information $D_{k-1}$ has 104 unequally spaced values from 1 to 1000. The stage cost $C_{S}$ has 35 unequally spaced values from 1 to 1000. The approximately optimal $n_{k}$ is chosen among of 368 candidate values $2,4,6,\ldots,180000,190000,200000$. The obtained data contain the approximately optimal $n_{k}$ for 3640 $(D_{k-1},C_{S})$ pairs. Figure~\ref{fig:contourcloglog} presents the contours of $n_{k}$ in the $(D_{k-1},C_{S})$-plane for the cloglog model. Figure~\ref{fig:contourcloglog} suggests that the contours are nearly linear in the log-scale when $D_{k-1}>10$.

The following model is found to have an excellent fit with the data
\begin{equation} \label{eq:costmodel}
\log n_{k}=\alpha+\beta \log D_{k-1}+\gamma \log C_{S}+\delta \log D_{k-1} \log C_{S}.
\end{equation}
The estimated model parameters $\alpha$, $\beta$, $\gamma$ and $\delta$ for the logit, the probit and the cloglog model are presented in Table~\ref{tab:models}. Apart from the intercept term, the estimated models are quite similar.

\begin{table}[htb]
\caption{Summary of estimated models~\eref{eq:costmodel}.  \label{tab:models} }
\begin{center}
\linespread{1}
\selectfont
\begin{tabular}{lcc}
\textbf{logit} & \multicolumn{2}{c}{$R^{2}=0.9992$}\\
& estimate & standard error\\
$\alpha$                      &   1.01515  &  0.00553 \\
$\beta$                  &   1.43396  &  0.00162 \\
$\gamma$                   &   0.41042  &  0.00134 \\
$\delta$      &   -0.01388  &  0.00039 \\
\\
\textbf{probit}  & \multicolumn{2}{c}{$R^{2}=0.9990$}\\
& estimate & standard error\\
$\alpha$                      &  0.12859  &  0.00611 \\
$\beta$                  &  1.41891  &  0.00179 \\
$\gamma$                   &  0.40324  &  0.00148 \\
$\delta$      &  -0.00949  &  0.00043 \\
\\
\textbf{cloglog} & \multicolumn{2}{c}{$R^{2}=0.9934$}\\
& estimate & standard error\\
$\alpha$                      &   0.48044  &  0.01366 \\
$\beta$                  &   1.34593  &  0.0042 \\
$\gamma$                   &   0.39711  &  0.00331 \\
$\delta$      &   -0.00778  &  0.00102
\end{tabular}
\end{center}
\end{table}

In the model~\eref{eq:costmodel} it is assumed that $a=1$. When $a$ is estimated, model~\eref{eq:costmodel} is written as
\begin{equation} \label{eq:costmodel2}
\log n_{k}=\alpha+\beta \log (\hat{a}D_{k-1})+\gamma \log C_{S}+\delta \log (\hat{a}D_{k-1}) \log C_{S}.
\end{equation}
In practical experiments, the stage size $n_{k}$ calculated from \eref{eq:costmodel2} need to be rounded to the closest even integer.

\section{Application to switching measurements} \label{sec:application}
\begin{figure}
  \includegraphics[width=\columnwidth]{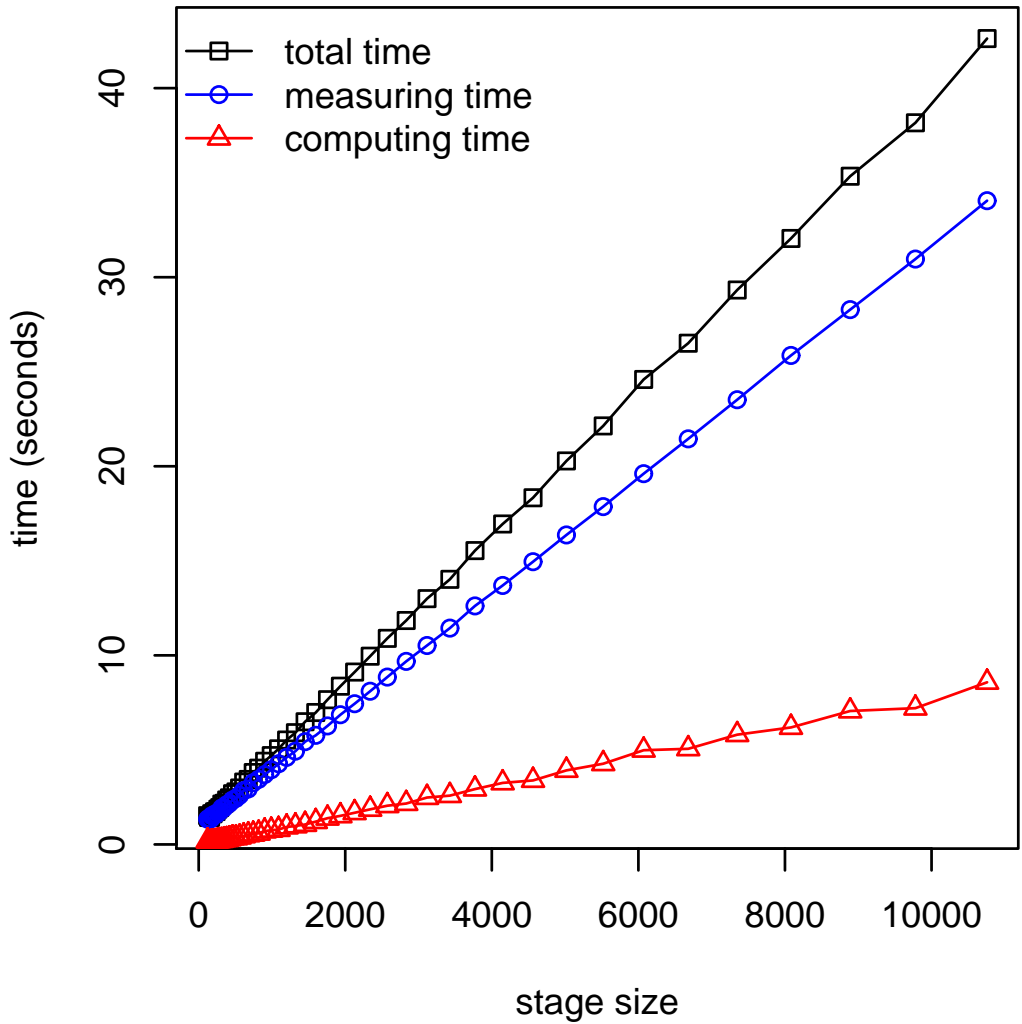}
  \caption{Time needed for measuring and for computing as function of stage size in the studied switching measurement.} \label{fig:times}
\end{figure}

The measurement data from \citep{optdesign} is reanalyzed in order to illustrate the efficiency of the presented approach. The data are available from the Royal Statistical Society Datasets Website \url{http://www.blackwellpublishing.com/rss/Volumes/Cv56p2.htm}. In the experiment, a sample consisting of aluminium--aluminium oxide--aluminium Josephson junction  circuit in a dilution refrigerator at 20~millikelvin temperature was connected to computer controlled measurement electronics in order to apply the current pulses and record the resulting voltage pulses. The initial maximum likelihood estimates were found using a binary search algorithm \citep{optdesign,inidesign}. After that the experiment continued as sequential design where D-optimal covariates values were calculated from the current maximum likelihood estimates of the parameters $a$ and $b$.  The approximate cost-efficient solution was not available at the time when the experiment was carried out and the stage size was set according to an ad-hoc rule where the stage size was 100 at the beginning and was then increased by 10 percent at each stage. The recorded data contains parameter estimates after each stage. In addition, the data contains the times needed for measuring and for computation. These times are plotted in Figure~\ref{fig:times} as a function of stage size. It can be seen that there is a linear relationship with the stage size and the total time. This supports the usage of the cost model defined in Section~\ref{sec:problemdefinition}. From the data it is estimate that the stage cost is 0.88167 seconds and the marginal cost of making one additional measurement is 0.00386 seconds. These numbers depend on the technical details of the measurement devices and cannot be generalized to other experiments. In order to apply the model~\eref{eq:costmodel2} the stage cost need to be standardized $C_{S}=0.88167/0.00386=228.4$.

\begin{figure}
  \includegraphics[width=\columnwidth]{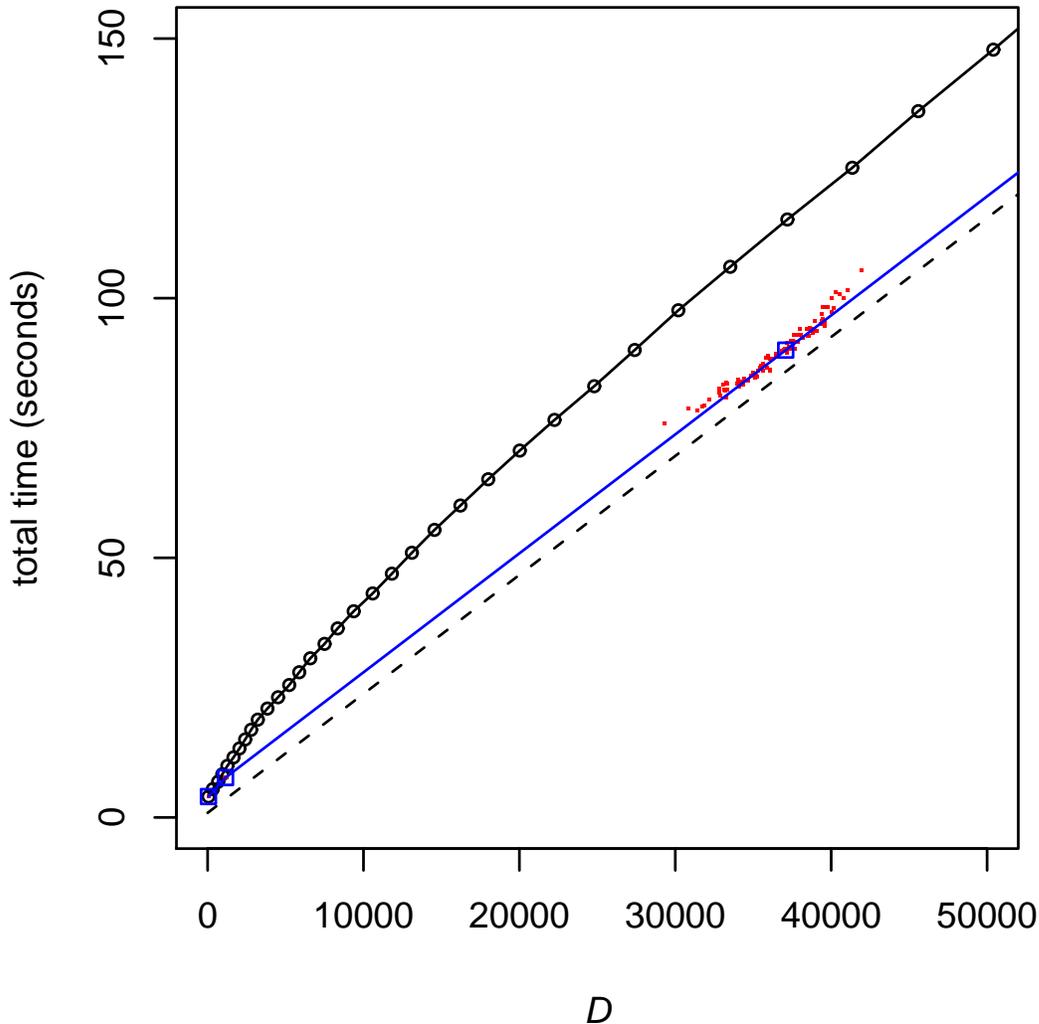}
  \caption{Comparison of the used sequential design and the cost-efficient sequential design in the $(D,C)$-plane. The circles show the stages of the design used in the actual switching measurement experiment. The squares show the path with the cost-efficient design calculated as the median of 100 simulated paths. The dots present the stages of the 100 simulated paths and give an idea on the variation between the simulation runs. It is instantly seen that the cost-efficient design saves time needed for the experiment. The dashed line shows the theoretical limit where the optimal covariate values are known from the beginning of the experiment.} \label{fig:switchingpaths}
\end{figure}

The design used in the experiment is compared to the cost-efficient design where the stage size is determined according to model~\eref{eq:costmodel}. The initial estimation took 4.0 seconds and resulted estimates $\hat{a}_{0}=0.380$, $\hat{b}=-95.60$ and $D_{0}=54.05$. Now model~\eref{eq:costmodel2} suggest the stage size 720 for the next stage while the stage size used in the experiment was only 134. Thus, the used ad-hoc rule leads to some losses in cost-efficiency and should not be used when it is important to minimize the time needed for the experiment. Further comparisons are carried out by means of simulation. For the simulation, it is assumed that the final estimates $\hat{a}=0.240$ and $\hat{b}=-60.628$ are the true parameter values. Starting from the initial estimates, the D-optimal covariate values and the cost-efficient stage size are calculated and new measurement data are generated. The simulation continues sequentially until $D$ exceeds the final $D$ achieved in the experiment. The resulting design path in $(D,C)$-plane calculated as the median of 100 simulated paths are presented in Figure~\ref{fig:switchingpaths}. It can be seen that the cost-efficient design is close to theoretical design where the D-optimal covariate values are known from the beginning and only one stage is needed. It is also seen that the ad-hoc rule leads to updating the parameter estimates too frequently. The total time needed for the experiment was 497 seconds, of which the cost-efficient design could have saved about 40 seconds.  This is not a crucial saving in a single switching measurement experiment but has importance when a number of experiments are carried out.

\section{Discussion} \label{sec:discussion}
Motivated by a practical problem from experimental physics, a framework for cost optimization in sequential designs has been formulated. Particularly, the approximately cost-efficient sequential designs for the logit, the probit and the cloglog model are determined. Cost-efficient designs for other link functions can be obtained in a similar way. The results can be directly applied to switching measurements and they may lead to substantial savings in the time needed for the experiment. In addition, the results are potentially useful in other experiments in physics, chemistry or engineering. Finding the specific applications will be a part of the future research.

The proposed approach may be criticized for using several approximations: Information matrix $\mathbf{J}$ is compressed to its determinant $D$. The change in $D$ is approximated by its expectation which itself is calculated numerically. The choice between design paths is based on a two-stage approximation. Finally, the optimal $n_{k}$ is obtained from the statistical model estimated from numerical data relying on all previous approximations. Naturally, the approximations have been used because analytical results are difficult to obtain. Nevertheless, it seems that the analytical results would not significantly improve cost-efficiency compared to the approximate results. Simulations in Section~\ref{sec:application} suggest that the approximate optimal design performed well compared to the theoretical D-optimal design with known parameters.

The idea of balancing by \citet{Sitter:optimaltwostage} and \citet{Sitter:twostagequantal} could be worth trying at least in the early stages of a sequential experiment. The calculation of the balanced design complicates the approach and takes some computational time but might result in further savings in the total time.

The concept of cost-efficiency narrows the gap between the theory of optimal design and the requirements of practical applications. When the cost is minimized, the most relevant criterion for the experiment is optimized. The presented cost model with a fixed stage cost can be easily extended to serve the needs of different applications. Due to the complexity of cost optimization, the determination of cost-efficient designs requires numerical approximations.  The numerical calculations may be complicated, but as it was demonstrated, once the optimization is done, the results are easily applied to all experiments with similar cost structure.

\section*{Ackowledgement}
The author thanks Dr. Juha J. Vartiainen and Dr. Teemu Pennanen for useful comments.

\end{document}